\documentclass{opt2024} 

\DeclareMathOperator*{\argmax}{arg\,max}

\title[Dimensionality Reduction Techniques for Global Bayesian Optimisation]{Dimensionality Reduction Techniques for Global Bayesian Optimisation}

\optauthor{%
\Name{Luo Long} \Email{luo.long@maths.ox.ac.uk}\\
\Name{Coralia Cartis} \Email{cartis@maths.ox.ac.uk}\\
\Name{Paz Fink Shustin} \Email{paz.finkshustin@maths.ox.ac.uk}\\
\addr Mathematical Institute, University of Oxford, Oxford, UK}

\begin{document}

\maketitle

\begin{abstract}%
Bayesian Optimisation (BO) is a state-of-the-art global optimisation technique for black-box problems where derivative information is unavailable and sample efficiency is crucial. However, improving the general scalability of BO has proved challenging. Here, we explore Latent Space Bayesian Optimisation (LSBO), that applies dimensionality reduction to perform BO in a reduced-dimensional subspace. While early LSBO methods used (linear) random projections (Wang et al., 2013 \cite{Wang2013}), we employ Variational Autoencoders (VAEs) to manage more complex data structures and general DR tasks. Building on Grosnit \textit{et al.} (2021) \cite{grosnit2021highdimensionalbayesianoptimisationvariational}, we analyse the VAE-based LSBO framework, focusing on VAE retraining and deep metric loss. We suggest a few key corrections in their implementation, originally designed for tasks such as molecule generation, and reformulate the algorithm for broader optimisation purposes. Our numerical results show that structured latent manifolds improve BO performance. Additionally, we examine the use of the Mat\'{e}rn-$\frac{5}{2}$ kernel for Gaussian Processes in this LSBO context. We also integrate Sequential Domain Reduction (SDR), a standard global optimization efficiency strategy, into BO. SDR is included in a GPU-based environment using \textit{BoTorch}, both in the original and VAE-generated latent spaces, marking the first application of SDR within LSBO.

\end{abstract}

\section{Introduction}
Global Optimisation (GO) aims to find the (approximate) global optimum of a smooth function \( f \) within a region of interest, possibly without the use of derivative problem information and with careful handling of often-costly objective evaluations.
In particular, we focus on the GO problem,
\begin{equation}\label{main problem}
    f^{\ast} = \min_{\mathbf{x} \in \mathcal{X}} f(\mathbf{x}), \tag{P}
\end{equation}
where \( \mathcal{X} \subseteq \mathbb{R}^D \) represents a feasible region, and \( f \) is a black-box, continuous function in (high) dimensions \( D \).
Bayesian Optimisation (BO) is a state-of-the-art GO framework that constructs a probabilistic model, typically a Gaussian Process (GP), of \( f \), and uses an acquisition function to guide sampling and efficiently search for the global optimum \cite{Frazier2018}. BO balances exploration and exploitation but suffers from scalability issues in high-dimensions \cite{hvarfner2024vanillabayesianoptimizationperforms}. To mitigate this, Dimensionality Reduction (DR) techniques can be used, allowing BO to operate in a lower-dimensional subspace where it is more effective. 
\paragraph{Contributions.} We investigate scaling up BO using DR techniques, focusing on VAEs and the Random Embedding Global Optimisation (REGO) framework \cite{cartis2021globaloptimizationusingrandom}. Our work extends the algorithm \cite{grosnit2021highdimensionalbayesianoptimisationvariational}, incorporating the Matérn-$\frac{5}{2}$ kernel, which enhances the flexibility and robustness of the method. Moreover, we conduct a comparative analysis of these approaches with standard BO techniques enhanced by Sequential Domain Reduction (SDR) \cite{SDR}. The first contribution of this work is that we \textbf{propose and implement SDR with BO within the \textit{BoTorch} framework \cite{balandat2020botorch}}, utilising GPU-based computation for efficiency. Furthermore, we \textbf{propose three BO-VAE algorithms, two of which are innovatively combined with SDR in the VAE-generated latent space} to boost optimisation performance. This marks the first instance of SDR being integrated with BO in the context of VAEs. We also investigate the effects of VAE retraining \cite{Tripp2020} and deep metric loss \cite{grosnit2021highdimensionalbayesianoptimisationvariational} on the optimisation process, emphasising the advantages of having a well-structured latent space for improved performance. Finally, we \textbf{compare our BO-VAE algorithms with the REMBO method} \cite{Wang2013} on low effective dimensionality problems, evaluating VAEs versus random embeddings as two different DR techniques in terms of optimisation performance.

\section{Preliminaries}
\paragraph{Bayesian Optimisation.} BO relies on two fundamental components: a GP prior and an acquisition function. Given a dataset of size $n$, $\mathcal{D}_n = \{\mathbf{x}_i, f(\mathbf{x}_i)\}_{i=1}^n$, the function values $\mathbf{f}_{1:n}$ are modelled as realisations of a Gaussian Random Vector (GRV) $\mathbf{F}_{1:n}$ under the GP prior. The distribution is characterised by a mean $\mathbb{E}_{\mathbf{F}_{1:n}}$ and covariance $K_{\mathbf{F}_{1:n}\mathbf{F}_{1:n}}$, where $\mathbf{F}_{1:n} \sim \mathcal{N}(\mathbb{E}_{\mathbf{F}_{1:n}}, K_{\mathbf{F}_{1:n}\mathbf{F}_{1:n}})$, using the Matérn-5/2 kernel. For an arbitrary unsampled point $\mathbf{x}$, the predicted function value $f(\mathbf{x})$ is inferred from the posterior distribution:
$
F(\mathbf{x}) \sim \mathcal{N}(\mu(\mathbf{x}|\mathcal{D}_n), \sigma^2(\mathbf{x}|\mathcal{D}_n)),
$
with
$
\mu(\mathbf{x}|\mathcal{D}_n) = \mathbb{E}_{F} + K_{\mathbf{F}_{1:n}F}^T K_{\mathbf{F}_{1:n}\mathbf{F}_{1:n}}^{-1}(\mathbf{f}_{1:n} - \mathbb{E}_{\mathbf{F}_{1:n}}),
$
$\sigma^2(\mathbf{x}|\mathcal{D}_n) = K_{FF} - K_{\mathbf{F}_{1:n}F}^T K_{\mathbf{F}_{1:n}\mathbf{F}_{1:n}}^{-1}K_{\mathbf{F}_{1:n}F}.$ BO uses the posterior mean $\mu(\cdot)$ and variance $\sigma^2(\cdot)$ in an acquisition function to guide sampling. In this work, we focus on Expected Improvement (EI):
$
u(\mathbf{x}|\mathcal{D}_n) = \mathbb{E}[\max\{F(\mathbf{x}) - f_n^{\ast}, 0\}|\mathcal{D}_n],
$
where $f_n^{\ast} = \max_{m \leq n} f(\mathbf{x}_m)$ is the highest observed value. To enhance BO, we incorporate SDR to refine the search region based on the algorithm’s best-found values, updating the region every few iterations to avoid missing the global optimum. To accelerate BO process, we propose to implement SDR \cite{SDR} within the traditional BO framework such that the search region can be refined to locate the global minimiser more efficiently according to the minimum function values found so far by the algorithm. Compared to the traditional SDR implementation that updates the search region at each iteration, we propose updating the region after a set number of iterations to avoid premature exclusion of the global optimum. Algorithm \ref{BOSDR} in Appendix \ref{Appendix BO_BOVAE with SDR} outlines the BO-SDR approach.

\paragraph{Variational Autoencoders.} DR methods reduce the number of features in a dataset while preserving essential information \cite{Velliangiri2019}. They can often be framed as an \textit{Encoder-Decoder} process, where the \textit{encoder} maps high-dimensional (HD) data to a lower-dimensional latent space, and the \textit{decoder} reconstructs the original data. We focus on VAEs \cite{kingma2022autoencodingvariationalbayes, doersch2021tutorialvariationalautoencoders}, a DR technique using Bayesian Variational Inference (VI) \cite{Hinton1993, Jordan1998}. VAEs utilise neural networks as encoders and decoders to generate latent manifolds. The probabilistic framework of a VAE consists of the encoder $q_{\boldsymbol{\phi}}(\cdot|\mathbf{x}): \mathcal{X} \rightarrow \mathcal{Z}$ parameterised by $\boldsymbol{\phi}$ which turns an input data $\mathbf{x} \in \mathbb{R}^D$ from some distribution into a distribution on the latent variable $\mathbf{z} \in \mathbb{R}^d$ ($d \ll D$), and the decoder $p_{\boldsymbol{\theta}}(\cdot|\mathbf{z}): \mathcal{Z} \rightarrow \mathcal{X}$ parameterised by $\boldsymbol{\theta}$ which reconstructs $\mathbf{x}$ as $\hat{\mathbf{x}}$ given samples from the latent distribution. The VAE's objective is to maximise the Evidence Lower BOund (ELBO):$ \mathcal{L}(\mathbf{\boldsymbol{\theta}, \boldsymbol{\phi}}; \mathbf{x})  =   \ln p_{\boldsymbol{\theta}}(\mathbf{x}) - D_{KL} [q_{\boldsymbol{\phi}}(\mathbf{z|x}) \| p_{\boldsymbol{\theta}}(\mathbf{z|x})] = \mathbb{E}_{q_{\boldsymbol{\phi}}(\mathbf{z|x})} [\ln p_{\boldsymbol{\theta}}(\mathbf{x|z})] - D_{KL} [q_{\boldsymbol{\phi}}(\mathbf{z|x}) \| p(\mathbf{z})],$ where $\ln p_{\boldsymbol{\theta}}(\mathbf{x})$ is the marginal log-likelihood, and $D_{KL} (\cdot \| \cdot)$ is the non-negative Kullback-Leibler Divergence between the true and the approximate posteriors. The prior is usually set to $\mathcal{N}(\textbf{0}, \mathbf{I})$, and the posterior is parametrised as Gaussians with diagonal covariance matrices, making ELBO optimisation tractable via the "reparameterisation trick" \cite{kingma2022autoencodingvariationalbayes}. Given $q_{\boldsymbol{\phi}}(\mathbf{z}|\mathbf{x}) = \mathcal{N}(\boldsymbol{\mu}(\mathbf{x}), \boldsymbol{\Sigma}(\mathbf{x}))$, the latent variable $\mathbf{z}$ is sampled as $\mathbf{z} = \mathbf{\boldsymbol{\mu}(x)} + \mathbf{\Sigma(x)}\boldsymbol{\xi}, \ \boldsymbol{\xi} \sim \mathcal{N}(\textbf{0}, \mathbf{I})$, enabling gradient-based optimisation with Adam \cite{kingma2017adammethodstochasticoptimization}.


\section{Algorithms}
As mentioned above, DR techniques help reduce the optimisation problem's dimensionality. Using a VAE within BO allows standard BO approach to be applied to larger scale problems, as then, we solve a GP regression sub-problem in the generated (smaller dimensional) latent space $\mathcal{Z}$. The BO-VAE approach\footnote{For brevity, we use BO-VAE to refer the approach of combing VAEs with BO.}, instead of solving \eqref{main problem} directly, attempts to solve
\begin{equation}\label{BO-VAE objective}
    f^\ast = \min_{\mathbf{z}\in \mathcal{Z}} \mathbb{E}_{p_{\boldsymbol{\theta}^\ast}\mathbf{(x|z)}} \left[f(\mathbf{x})\right], 
\end{equation}
where $\boldsymbol{\theta}^\ast$ is the optimal decoder network parameter. Therefore, it is implicitly assumed that the optimal point $\mathbf{x}^\ast$ can be obtained from the optimal decoder with some probability given some latent data $\mathbf{z}$ by the associated optimal encoder $q_{\boldsymbol{\phi}^\ast}\mathbf{(z|x)}$ \cite{grosnit2021highdimensionalbayesianoptimisationvariational}, $ \exists \ \mathbf{z}\in\mathcal{Z}, \mathbb{P}\left[ \mathbf{x}^\ast \sim p_{\boldsymbol{\theta}^\ast}\mathbf{(\cdot|z)}\right] > 0.$ 

When fitting the GP surrogate, we follow \cite{grosnit2021highdimensionalbayesianoptimisationvariational} and use Deep Metric Loss (DML) to generate well-structured VAE-generated latent spaces. Specifically, we apply the soft triplet loss and retrain the VAEs following \cite{Tripp2020} to adapt to new points from the GP and optimise the black-box objective efficiently. Additionally, we implement SDR in the latent space to accelerate the BO process. Algorithm \ref{general BOVAE} outlines our BO-VAE approach with SDR, consisting of the pre-training of a standard VAE on the unlabelled dataset $\mathcal{D}_{\mathbb{U}}$ (line 1) and optional retraining with soft triplet loss to structure the latent space by gradually adjusting the network parameters of the encoder and decoder. When the soft triplet loss is used in retraining the VAE, the modified VAE ELBO $\mathcal{L}_{DML}(\cdot)$ is used in line 4 instead; see Appendix \ref{Appendix triplet loss} for details. The BO-VAE algorithm with DML is included in Appendix \ref{Appendix triplet loss} as Algorithm \ref{BOVAE with DML}. For comparison, we provide a baseline BO-VAE algorithm without retraining or DML (Appendix \ref{Appendix BO_BOVAE with SDR}, Algorithm \ref{BOVAE-SDR}). Theorem 1 in \cite{grosnit2021highdimensionalbayesianoptimisationvariational} offers a regret analysis with a sub-linear convergence rate, providing a valuable theoretical foundation. However, the proof relies on the assumption of a Gaussian kernel, limiting its direct applicability when using the Matérn kernel, as we do here. Despite this limitation, the theorem provides key insights supporting the BO-VAE approach. Our ongoing work addresses this gap, and a similar result specifically tailored to the Matérn kernel is delegated to future work.
\begin{algorithm2e}[!htb]\label{general BOVAE}
  \small
  \caption{Retraining BO-VAE Algorithm with SDR}
  \SetAlgoLined
  \LinesNumbered
  \KwData{Labelled dataset $\mathcal{D}_\mathbb{L}^{l = 1} = \{\mathbf{x}_i, f(\mathbf{x}_i)\}_{i = 1}^N$, unlabelled dataset $\mathcal{D}_\mathbb{U} = \{\mathbf{x}_i\}_{i = 1}^M$, budget $B$, periodic frequency $q$, initial bound $R^0$ in latent space $\mathcal{Z}$, EI acquisition function $u(\cdot)$, the encoder and decoder models from a VAE, $q_{\boldsymbol{\phi}}\mathbf{(z|x)}: \mathcal{X} \rightarrow \mathcal{Z}$ and $p_{\boldsymbol{\theta}}\mathbf{(x|z)}: \mathcal{Z} \rightarrow \mathcal{X}$.}
  \KwResult{Minimum function value $f_{min}$ found by the algorithm.}

  \textbf{Pre-train} the VAE model $V_{\mathcal{D}_\mathbb{L}}^{l=0}$ with $\mathcal{D}_\mathbb{U}$: $
  \boldsymbol{\theta}^\ast_0, \boldsymbol{\phi}^\ast_0 = \argmax_{\boldsymbol{\theta}, \boldsymbol{\phi}} \mathcal{L}(\boldsymbol{\theta}, \boldsymbol{\phi}; \mathcal{D}_\mathbb{U})
  $\;
  
  Set $\boldsymbol{\theta}^\ast_1 \leftarrow  \boldsymbol{\theta}^\ast_0$, $\boldsymbol{\phi}^\ast_1 \leftarrow \boldsymbol{\phi}^\ast_0$, $V_{\mathcal{D}_\mathbb{L}}^{l=1} \leftarrow V_{\mathcal{D}_\mathbb{L}}^{l=0}$\;

  \For{$l = 1$ \text{to} $L \equiv \lceil B/q \rceil$}{
    \textbf{Train} the VAE model $V_{\mathcal{D}_\mathbb{L}^l}^{l}$ on $\mathcal{D}_\mathbb{L}$:
    $
    \boldsymbol{\theta}^\ast_l, \boldsymbol{\phi}^\ast_l = \argmax_{\boldsymbol{\theta}, \boldsymbol{\phi}} \mathcal{L}(\boldsymbol{\theta}, \boldsymbol{\phi}; \mathcal{D}_\mathbb{L}^l)
    $\;

    \textbf{Compute} the latent dataset 
    $
    \mathcal{D}_\mathbb{Z}^l = \{ \mathbf{z}_i, f(\mathbf{x}_i) \}^{N + l \cdot q}_{i = 1} = \{ \mathbb{E}_{q_{\boldsymbol{\phi}^\ast_l}(\mathbf{z|x_i})}[\mathbf{z}], f(\mathbf{x}_i) \}^{N + l \cdot q}_{i = 1}
    $\;
    
    \textbf{Initialise} $\mathcal{D}_\mathbb{L}^{l;k=0} \leftarrow \mathcal{D}_\mathbb{L}^{l}$ and $\mathcal{D}_\mathbb{Z}^{l;k=0} \leftarrow \mathcal{D}_\mathbb{Z}^{l}$\;
    \textbf{Initialise} SDR with $R^0$\;
    \For{$k = 0$ \text{to} $q-1$}{
      \textbf{Fit} a Gaussian Process (GP) model $h_{l;k}: \mathcal{Z} \rightarrow \mathbb{R}$ on
      $
      \mathcal{D}_\mathbb{Z}^{l;k} = \{ \mathbf{z}_i, f(\mathbf{x}_i) \}^{N + l \cdot q + k}_{i = 0}
      $\;
      
      \textbf{Solve} for the next latent point:
      $
      \hat{\mathbf{z}}_{l; k + 1} = \argmax_{\mathbf{z}} u(\mathbf{z}|\mathcal{D}_\mathbb{Z}^{l;k})
      $\;
      
      \textbf{Obtain} the new sample $\hat{\mathbf{x}}_{l; k + 1}$: 
      $
      \hat{\mathbf{x}}_{l; k + 1} \sim p_{\boldsymbol{\theta}^\ast_l}(\cdot| \hat{\mathbf{z}}_{l; k + 1})
      $\;

      \textbf{Evaluate} the objective function at the new sample: $f(\hat{\mathbf{x}}_{l; k + 1})$\;

      \textbf{Augment} the datasets:
      $
      \mathcal{D}_\mathbb{L}^{l; k+1 } \leftarrow \mathcal{D}_\mathbb{L}^{l;k} \cup \{ \hat{\mathbf{x}}_{l; k+ 1}, f(\hat{\mathbf{x}}_{l; k+1}) \}
      , \ 
      \mathcal{D}_\mathbb{Z}^{l; k+1} \leftarrow \mathcal{D}_\mathbb{Z}^{l;k} \cup \{ \hat{\mathbf{z}}_{l; k + 1}, f(\hat{\mathbf{x}}_{l; k + 1}) \}$\;
      
      \textbf{Update} the search domain:
      $
      R^{k + 1} \leftarrow R^k
      $
      using SDR given $\mathcal{D}_{\mathbb{Z}}^{l; k+1}$ \;
    }

    \textbf{Augment} the outer-loop datasets:
    $
    \mathcal{D}_\mathbb{L}^{l+1} \leftarrow \mathcal{D}_\mathbb{L}^{l;q}, \mathcal{D}_\mathbb{Z}^{l+1} \leftarrow \mathcal{D}_\mathbb{Z}^{l;q}
    $\;
  }
\end{algorithm2e}
\section{Numerical Study}

We conduct numerical experiments with the three BO-VAE algorithms (Algorithms \ref{general BOVAE}, \ref{BOVAE-SDR}, \ref{BOVAE with DML}) within the \textit{BoTorch} framework \cite{balandat2020botorch}. We explore cases where $d = 2, 5$ for $D = 10$, and $d = 2, 10, 50$ for $D = 100$. The results reveal that, for a fixed ambient dimension $D$, larger latent dimensions, particularly $d = 50$, tend to degrade performance due to the reduced VAE generalisation capacity. In contrast, smaller latent dimensions ($d = 2, 5$) yield better results, as VAEs then can give more efficient latent data representations, and the BO can solve such reduced problems more efficiently. In this paper, we present experimental results for the case $D = 100$ and $d = 2$, which strongly highlight the advantages of SDR in latent space optimisation and illustrate the three BO-VAE algorithms. The encoder structure of the VAE used is $[100, 30, 2]$\footnote{It indicates a three-layer feedforward neural network: the input layer has $100$ neurons, followed by a hidden layer with $30$ neurons, and finally an output layer with $2$ neurons. Similarly for the others.}, and the decoder structure is $[2, 30, 100].$ For the REMBO comparison, the VAE used has $[100, 25, 5]$ for the encoder and $[5, 25, 100]$ for the decoder. The activation function is Soft-plus. We set the budget $B = 350$ and $q = 50$ to retrain 7 times for Algorithms \ref{general BOVAE} and \ref{BOVAE with DML}. In practice, to improve computational costs, the input spaces of VAEs can be normalised to a fixed range, such as a hypercube, simplifying both pre-training and retraining steps and avoiding the need to train multiple VAEs. By applying an additional scaling between the specific problem domain and the fixed VAE input space, we leverage the VAE’s generalisation ability to construct latent manifolds for datasets from diverse problem domains. In our experiments, we set the fixed VAE input space to $[-3, 3]^D$. Further details of our experiments are provided in Appendix \ref{Appendix Additional Details for Section 4 Numerical Experiments}.

\begin{figure}[!htb]
    \centering
    \includegraphics[width=0.67\linewidth]{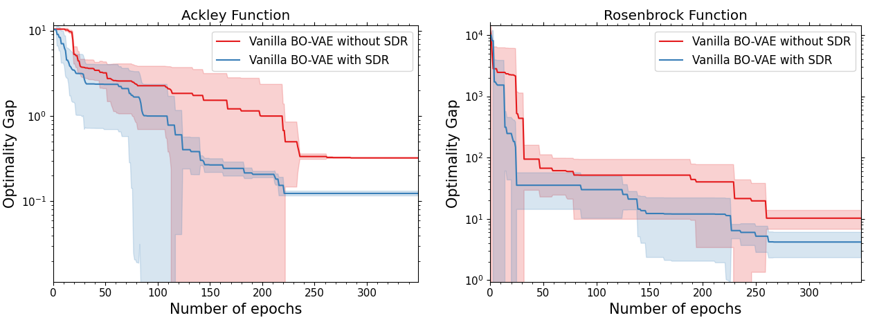}
    \caption{\small Comparisons of Vanilla BO-VAE algorithm (Algorithm \ref{BOVAE-SDR}) with and without SDR on $100$D Ackley and Rosenbrock problems. The means and the standard deviations (shaded areas) of the minimum function values found are plotted across $5$ repeated runs.}
    \label{fig: SDR in VAE}
\end{figure}

\begin{figure}[!htb]
    \centering
    \includegraphics[width=0.67\linewidth]{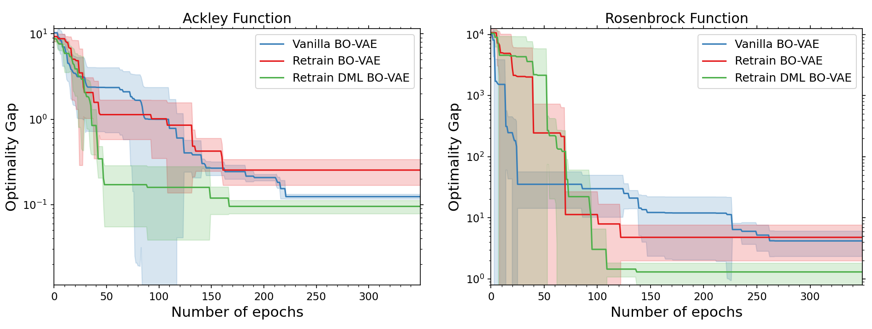}
    \caption{\small Comparisons of Algorithm \ref{BOVAE-SDR} (Vanilla BO-VAE), Algorithm \ref{general BOVAE} (Retrain BO-VAE), and Algorithm \ref{BOVAE with DML} (Retrain DML BO-VAE) in solving $100$D Ackley and Rosenbrock problems. The means and the standard deviations (shaded areas) of the minimum function values found are plotted across $5$ repeated runs.}
    \label{fig: D = 100}
\end{figure}

\paragraph{Numerical Illustration of SDR within BO-VAE.} To demonstrate the effectiveness of SDR within VAE-generated latent subspaces, we use Algorithm \ref{BOVAE-SDR} to generate results in Figure \ref{fig: SDR in VAE}. The figure shows that SDR leads to faster convergence and lower, potentially optimal, function values. Additionally, Figure \ref{fig: D = 100} illustrates the results of the three BO-VAE algorithms, where Algorithm \ref{BOVAE with DML} generally outperforms the others, achieving better global optima. This improvement is due to the soft triplet loss, which better structures the latent subspaces and enhances the GP surrogate's efficiency.
\paragraph{Algorithm Comparisons.} We first compare the three BO-VAE algorithms with BO-SDR on Test Set 1 (Table \ref{tab:Full-rank_Test_problem}) and then against REMBO \cite{Wang2013} on Test Set 2 (see Appendix \ref{Appendix High-dimensional Low-rank Test Set}). REMBO addresses \eqref{main problem} by solving a reduced problem in a low-dimensional subspace: $ \min_{\mathbf{y} \in \mathbb{R}^d} f(\mathbf{Ay}) = \min_{\mathbf{y} \in \mathbb{R}^d} g(\mathbf{y})$, subject to $ \mathbf{y} \in \mathcal{Y} = [-\delta, \delta]^d.$ Here $\mathbf{A}$ is a $D \times d$ Gaussian matrix for random embedding, with $d \ll D$. Solving $g(\mathbf{y})$ in the reduced subspace is equivalent to solving $f(\mathbf{Ay})$. In this context, the Gaussian matrix $\mathbf{A}$ serves as a (linear) encoder (for dimensionality reduction), while $\mathbf{A}^T$ acts as a decoder. For the REMBO comparisons, we used $d = d_e + 1$, where $d_e$ is the effective dimensionality, and set $\delta = 2.2\sqrt{d_e}$ based on \cite{cartis2020dimensionalityreductiontechniqueunconstrained}. 
Each (randomised) algorithm was run twice on each problem in Test Set 1. Each (randomised) problem in Test Set 2 was run twice, yielding 10 test problems. 
The algorithms we are comparing are denoted by BO-SDR (Algorithm \ref{BOSDR}), V-BOVAE (Algorithm \ref{BOVAE-SDR}), R-BOVAE (Algorithm \ref{general BOVAE}), and S-BOVAE (Algorithm \ref{BOVAE with DML}). 

\begin{figure}[!htb]
    \centering
    \includegraphics[scale = 0.30]{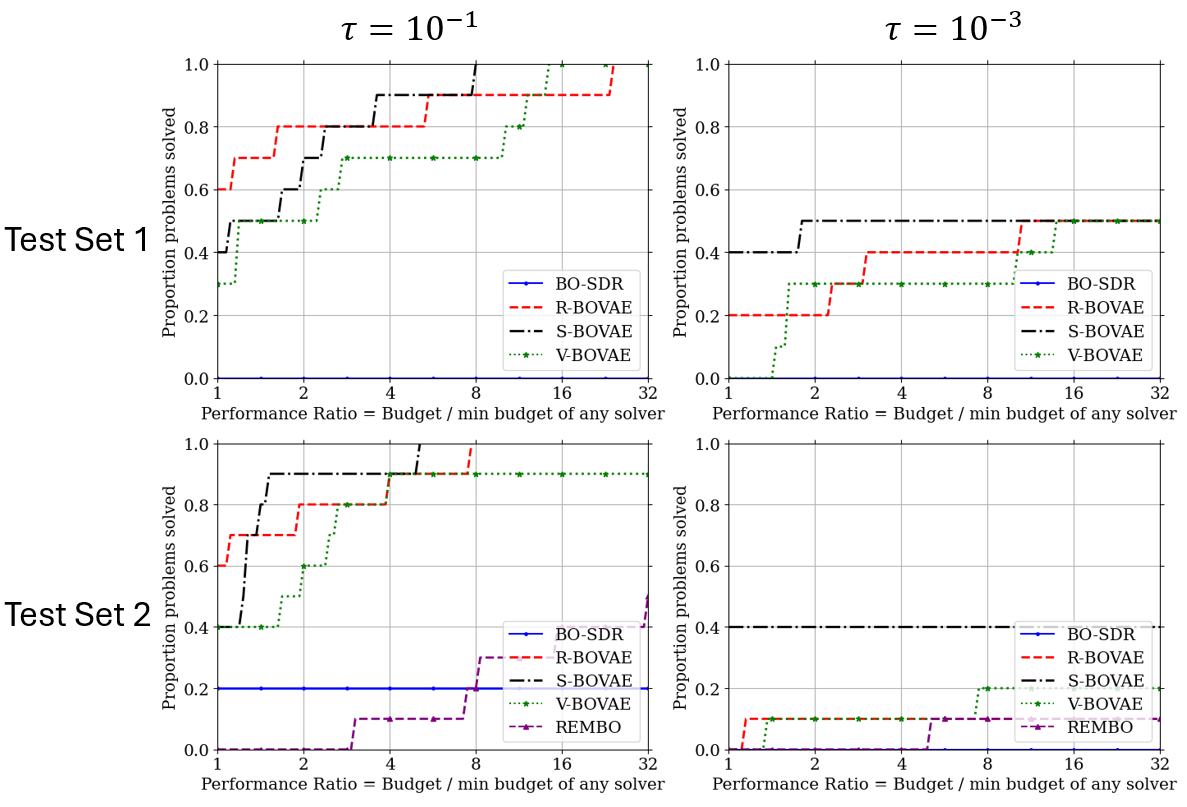}
    \caption{\small Performance Profiles on Test Sets 1 \& 2 when $\tau = 10^{-1}$ and $10^{-3}$.}
    \label{fig: Alg Comp}
\end{figure}

Results with accuracy levels \(\tau = 10^{-1}\) and \(\tau = 10^{-3}\) are summarised in Tables \ref{test_set1_data_profiles}, \ref{test_set2_data_profiles} and Figure \ref{fig: Alg Comp}. From these results, it can be  seen that BO-VAE algorithms solve more problems compared to BO-SDR and REMBO algorithms. While BO-SDR may struggle with scalability, REMBO's low problem-solving percentages are likely due to the over-exploration of the boundary projections \cite{AFrameworkforBayesianOptimizationinEmbeddedSubspaces} and the embedding subspaces failing to accurately capture the global minimisers. To address this, \cite{Wang2013} recommends restarting REMBO to improve the success rate. Meanwhile, Algorithm \ref{BOVAE with DML} (S-BOVAE) consistently performed best due to its structured latent spaces.
\begin{table}[!htb]
    \centering
    \begin{minipage}{.47\textwidth}
        \centering
        \resizebox{.60\textwidth}{!}{\begin{tabular}{| c |c  | c | }
        \hline
        & {$\tau = 10^{-1}$} & {$\tau = 10^{-3}$}\\
        \hline
        BO-SDR & $10\%$ & $0 \%$ \\
        \hline
        V-BOVAE & $100\%$ & $50\%$  \\
        \hline
        S-BOVAE & $100\%$  & $50\%$ \\
        \hline
        R-BOVAE & $100\%$ & $50\%$  \\
        \hline
        \end{tabular}}
        \caption{\small Average percentage of problems solved in Test Set 1 for $\tau = 10^{-1}$ and $10^{-3}$.}
        \label{test_set1_data_profiles}
    \end{minipage}
    \hfill
    \begin{minipage}{.47\textwidth}
        \centering
        \resizebox{.60\textwidth}{!}{\begin{tabular}{| c | c | c | c | c |}
        \hline
        & $\tau = 10^{-1}$ & $\tau = 10^{-3}$\\
        \hline
        BO-SDR & $20\%$   & $0 \%$  \\
        \hline
        V-BOVAE & $90\%$ & $20\%$   \\
        \hline
        S-BOVAE & $100\%$ & $40\%$   \\
        \hline
        R-BOVAE & $100\%$ & $30\%$  \\
        \hline
        REMBO & $50\%$ & $10\%$\\
        \hline
        \end{tabular}}
        \caption{\small Average percentage of problems solved in Test Set 2 for $\tau = 10^{-1}$ and $10^{-3}$.}
        \label{test_set2_data_profiles}
    \end{minipage}
\end{table}

\section{Conclusion and Future Work}
In this work, we have explored dimensionality reduction techniques to enhance the scalability of BO. The use of VAEs offers an alternative and more general approach for GP fitting in low-dimensional latent subspaces, alleviating the curse of dimensionality. Unlike REMBO, which primarily targets low-rank functions, VAE-based LSBO is effective for both high-dimensional full-rank and low-rank functions. Although BO-VAE reduces function values effectively, the optimality gap remains constrained by noise from the VAE loss, as seen in Figures \ref{fig: SDR in VAE} and \ref{fig: D = 100}. To address this, implementations of data weights and different GP initialisation are the potential future directions. Additionally, SDR struggles in high dimensions; adopting methods like domain refinement based on threshold probabilities \cite{salgia2021domainshrinkingbasedbayesianoptimization} may improve performance.\\

\paragraph{Acknowledgments} The second author's (CC) research was supported by  the Hong Kong Innovation and Technology Commission's Center for Intelligent Multi-dimensional Analysis (InnoHK Project CIMDA).

\bibliography{sample}

\begin{thebibliography}{28}
\providecommand{\natexlab}[1]{#1}
\providecommand{\url}[1]{\texttt{#1}}
\expandafter\ifx\csname urlstyle\endcsname\relax
  \providecommand{\doi}[1]{doi: #1}\else
  \providecommand{\doi}{doi: \begingroup \urlstyle{rm}\Url}\fi

\bibitem[Balandat et~al.(2020)Balandat, Karrer, Jiang, Daulton, Letham, Wilson, and Bakshy]{balandat2020botorch}
Maximilian Balandat, Brian Karrer, Daniel~R. Jiang, Samuel Daulton, Benjamin Letham, Andrew~Gordon Wilson, and Eytan Bakshy.
\newblock {BoTorch: A Framework for Efficient Monte-Carlo Bayesian Optimization}.
\newblock In \emph{Advances in Neural Information Processing Systems 33}, 2020.
\newblock URL \url{http://arxiv.org/abs/1910.06403}.

\bibitem[Bowman et~al.(2016)Bowman, Vilnis, Vinyals, Dai, Jozefowicz, and Bengio]{bowman2016generatingsentencescontinuousspace}
Samuel~R. Bowman, Luke Vilnis, Oriol Vinyals, Andrew~M. Dai, Rafal Jozefowicz, and Samy Bengio.
\newblock Generating sentences from a continuous space, 2016.
\newblock URL \url{https://arxiv.org/abs/1511.06349}.

\bibitem[Burgess et~al.(2018)Burgess, Higgins, Pal, Matthey, Watters, Desjardins, and Lerchner]{burgess2018understandingdisentanglingbetavae}
Christopher~P. Burgess, Irina Higgins, Arka Pal, Loic Matthey, Nick Watters, Guillaume Desjardins, and Alexander Lerchner.
\newblock Understanding disentangling in $\beta$-vae, 2018.
\newblock URL \url{https://arxiv.org/abs/1804.03599}.

\bibitem[Cartis and Otemissov(2020)]{cartis2020dimensionalityreductiontechniqueunconstrained}
Coralia Cartis and Adilet Otemissov.
\newblock A dimensionality reduction technique for unconstrained global optimization of functions with low effective dimensionality, 2020.
\newblock URL \url{https://arxiv.org/abs/2003.09673}.

\bibitem[Cartis et~al.(2021{\natexlab{a}})Cartis, Massart, and Otemissov]{cartis2021globaloptimizationusingrandom}
Coralia Cartis, Estelle Massart, and Adilet Otemissov.
\newblock Global optimization using random embeddings, 2021{\natexlab{a}}.
\newblock URL \url{https://arxiv.org/abs/2107.12102}.

\bibitem[Cartis et~al.(2021{\natexlab{b}})Cartis, Roberts, and Sheridan-Methven]{Cartis_2021}
Coralia Cartis, Lindon Roberts, and Oliver Sheridan-Methven.
\newblock Escaping local minima with local derivative-free methods: a numerical investigation.
\newblock \emph{Optimization}, 71\penalty0 (8):\penalty0 2343–2373, February 2021{\natexlab{b}}.
\newblock ISSN 1029-4945.
\newblock \doi{10.1080/02331934.2021.1883015}.
\newblock URL \url{http://dx.doi.org/10.1080/02331934.2021.1883015}.

\bibitem[Cartis et~al.(n.d.)Cartis, Fowkes, and Roberts]{CartisFowkesRoberts}
Coralia Cartis, Jaroslav Fowkes, and Lindon Roberts.
\newblock Optimization resources: A collection of software and resources for nonlinear optimization.
\newblock \url{https://lindonroberts.github.io/opt/resources.html#data-performance-profiles}, n.d.

\bibitem[Doersch(2021)]{doersch2021tutorialvariationalautoencoders}
Carl Doersch.
\newblock Tutorial on variational autoencoders, 2021.
\newblock URL \url{https://arxiv.org/abs/1606.05908}.

\bibitem[Ernesto and Diliman(2005)]{Ernesto2005}
P.A. Ernesto and U.P. Diliman.
\newblock {MVF}—multivariate test functions library in c for unconstrained global optimization, 2005.

\bibitem[Frazier(2018)]{Frazier2018}
Peter~I. Frazier.
\newblock A tutorial on bayesian optimization.
\newblock \emph{arXiv preprint arXiv:1807.02811}, 2018.
\newblock URL \url{https://arxiv.org/abs/1807.02811}.

\bibitem[Fu et~al.(2019)Fu, Li, Liu, Gao, Celikyilmaz, and Carin]{fu2019cyclicalannealingschedulesimple}
Hao Fu, Chunyuan Li, Xiaodong Liu, Jianfeng Gao, Asli Celikyilmaz, and Lawrence Carin.
\newblock Cyclical annealing schedule: A simple approach to mitigating kl vanishing, 2019.
\newblock URL \url{https://arxiv.org/abs/1903.10145}.

\bibitem[Grosnit et~al.(2021)Grosnit, Tutunov, Maraval, Griffiths, Cowen-Rivers, Yang, Zhu, Lyu, Chen, Wang, Peters, and Bou-Ammar]{grosnit2021highdimensionalbayesianoptimisationvariational}
Antoine Grosnit, Rasul Tutunov, Alexandre~Max Maraval, Ryan-Rhys Griffiths, Alexander~I. Cowen-Rivers, Lin Yang, Lin Zhu, Wenlong Lyu, Zhitang Chen, Jun Wang, Jan Peters, and Haitham Bou-Ammar.
\newblock High-dimensional bayesian optimisation with variational autoencoders and deep metric learning, 2021.
\newblock URL \url{https://arxiv.org/abs/2106.03609}.

\bibitem[Hinton and Camp(1993)]{Hinton1993}
G.~E. Hinton and D.~Van Camp.
\newblock Keeping the neural networks simple by minimizing the description length of the weights.
\newblock In \emph{Proceedings of the Sixth Annual Conference on Computational Learning Theory}, pages 5--13, 1993.

\bibitem[Hoffer and Ailon(2015)]{Hoffer2015}
Elad Hoffer and Nir Ailon.
\newblock Deep metric learning using triplet network.
\newblock In \emph{International Workshop on Similarity-Based Pattern Recognition}, pages 84--92. Springer, 2015.

\bibitem[Hvarfner et~al.(2024)Hvarfner, Hellsten, and Nardi]{hvarfner2024vanillabayesianoptimizationperforms}
Carl Hvarfner, Erik~Orm Hellsten, and Luigi Nardi.
\newblock Vanilla bayesian optimization performs great in high dimensions, 2024.
\newblock URL \url{https://arxiv.org/abs/2402.02229}.

\bibitem[Ishfaq et~al.(2018)Ishfaq, Hoogi, and Rubin]{Ishfaq2018}
Haque Ishfaq, Assaf Hoogi, and Daniel Rubin.
\newblock Tvae: Triplet-based variational autoencoder using metric learning.
\newblock \emph{arXiv preprint arXiv:1802.04403}, 2018.
\newblock URL \url{https://arxiv.org/abs/1802.04403}.

\bibitem[Jordan et~al.(1998)Jordan, Ghahramani, Jaakkola, and Saul]{Jordan1998}
M.~I. Jordan, Z.~Ghahramani, T.~S. Jaakkola, and L.~K. Saul.
\newblock An introduction to variational methods for graphical models.
\newblock In \emph{Learning in Graphical Models}, pages 105--161. Springer, 1998.

\bibitem[Kingma and Ba(2017)]{kingma2017adammethodstochasticoptimization}
Diederik~P. Kingma and Jimmy Ba.
\newblock Adam: A method for stochastic optimization, 2017.
\newblock URL \url{https://arxiv.org/abs/1412.6980}.

\bibitem[Kingma and Welling(2022)]{kingma2022autoencodingvariationalbayes}
Diederik~P Kingma and Max Welling.
\newblock Auto-encoding variational bayes, 2022.
\newblock URL \url{https://arxiv.org/abs/1312.6114}.

\bibitem[Moré and Wild(2009)]{More2009}
J.~J. Moré and S.~M. Wild.
\newblock Benchmarking derivative-free optimization algorithms.
\newblock \emph{SIAM Journal on Optimization}, 20:\penalty0 172--191, 2009.

\bibitem[Nayebi et~al.(2019)Nayebi, Munteanu, and Poloczek]{AFrameworkforBayesianOptimizationinEmbeddedSubspaces}
Amin Nayebi, Alexander Munteanu, and Matthias Poloczek.
\newblock A framework for bayesian optimization in embedded subspaces.
\newblock In Kamalika Chaudhuri and Ruslan Salakhutdinov, editors, \emph{Proceedings of the 36th International Conference on Machine Learning}, volume~97 of \emph{Proceedings of Machine Learning Research}, pages 4752--4761. PMLR, 09--15 Jun 2019.
\newblock URL \url{https://proceedings.mlr.press/v97/nayebi19a.html}.

\bibitem[Salgia et~al.(2021)Salgia, Vakili, and Zhao]{salgia2021domainshrinkingbasedbayesianoptimization}
Sudeep Salgia, Sattar Vakili, and Qing Zhao.
\newblock A domain-shrinking based bayesian optimization algorithm with order-optimal regret performance, 2021.
\newblock URL \url{https://arxiv.org/abs/2010.13997}.

\bibitem[Stander and Craig(2002)]{SDR}
Nielen Stander and Kenneth Craig.
\newblock On the robustness of a simple domain reduction scheme for simulation-based optimization.
\newblock \emph{International Journal for Computer-Aided Engineering and Software (Eng. Comput.)}, 19, 06 2002.
\newblock \doi{10.1108/02644400210430190}.

\bibitem[Surjanovic and Bingham(2013)]{Surjanovic2013}
S.~Surjanovic and D.~Bingham.
\newblock Virtual library of simulation experiments: Test functions and datasets.
\newblock \url{https://www.sfu.ca/~ssurjano/}, 2013.

\bibitem[Tripp et~al.(2020)Tripp, Daxberger, and Hernández-Lobato]{Tripp2020}
Austin Tripp, Erik Daxberger, and José~Miguel Hernández-Lobato.
\newblock Sample-efficient optimization in the latent space of deep generative models via weighted retraining.
\newblock In \emph{Advances in Neural Information Processing Systems}, volume~33, 2020.

\bibitem[Velliangiri et~al.(2019)Velliangiri, Alagumuthukrishnan, and Joseph]{Velliangiri2019}
S.~Velliangiri, S.~Alagumuthukrishnan, and S.~I.~Thankumar Joseph.
\newblock A review of dimensionality reduction techniques for efficient computation.
\newblock \emph{Procedia Computer Science}, 165:\penalty0 104--111, 2019.
\newblock \doi{10.1016/j.procs.2020.01.079}.
\newblock URL \url{https://doi.org/10.1016/j.procs.2020.01.079}.

\bibitem[Wang et~al.(2013)Wang, Zoghi, Hutter, Matheson, and Freitas]{Wang2013}
Z.~Wang, M.~Zoghi, F.~Hutter, D.~Matheson, and N.~De Freitas.
\newblock Bayesian optimization in high dimensions via random embeddings.
\newblock In \emph{International Joint Conference on Artificial Intelligence}, pages 1778--1784, 2013.

\bibitem[Wu et~al.(2019)Wu, Rallabandi, Black, and Nyberg]{Wu2019}
Peter Wu, SaiKrishna Rallabandi, Alan~W. Black, and Eric Nyberg.
\newblock Ordinal triplet loss: Investigating sleepiness detection from speech.
\newblock In \emph{Proc. Interspeech 2019}, pages 2403--2407, 2019.
\newblock \doi{10.21437/Interspeech.2019-2278}.
\newblock URL \url{http://dx.doi.org/10.21437/Interspeech.2019-2278}.

\end{thebibliography}
\newpage
\clearpage
\appendix

\section{Test Sets}



\subsection{High-dimensional Full-rank Test Set}

\begin{table}[!htb]
    \centering
    \begin{tabular}{|c|c|c|c|c|}
    \hline
    \# & Function & Dimension(s) & Domain & Global Minimum \\
    \hline
    1 & Ackley \cite{Ernesto2005} & $D$ & $\mathbf{x} \in [-30, 30]^{D}$ & 0 \\
    \hline 
    2 & Levy \cite{Surjanovic2013} & $D$ & $\mathbf{x} \in [-10, 10]^D$ & 0 \\
    \hline
    3 & Rosenbrock \cite{Surjanovic2013} & $D$ & $\mathbf{x} \in [-5, 10]^D$ & 0 \\
    \hline
    4 & Styblinski-Tang \cite{Surjanovic2013} & $D$ & $\mathbf{x} \in [-5, 5]^D$ & $-39.16599 \times D$ \\
    \hline
    5 & Rastrigin Function \cite{Surjanovic2013} & $D$ & $\mathbf{x} \in [-5.12, 5.12]^D$ & 0 \\
    \hline
    \end{tabular}
\caption{Benchmark high-dimensional full-rank test problems from \cite{Cartis_2021, CartisFowkesRoberts}.}
\label{tab:Full-rank_Test_problem}
\end{table}


\subsection{High-dimensional Low-rank Test Set}\label{Appendix High-dimensional Low-rank Test Set}
The low-rank test set, or Test Set 2, comprises $D$-dimensional low-rank functions generated from the low-rank test functions listed in Table \ref{Low-rank Test problem}. To construct these $D$-dimensional functions with low effective dimensionality, we adopt the methodology proposed in \cite{Wang2013}.\\
 Let $\Bar{h}(\Bar{\mathbf{x}})$ be any function from Table \ref{Low-rank Test problem} with dimension $d_e$ and the given domain scaled to $[-1, 1]^{d_e}$. The first step is to append $D -d_e$ fake dimensions with zero coefficients to $\Bar{h}(\Bar{\mathbf{x}})$:
\begin{equation*}
    h(\mathbf{x}) = \Bar{h}(\Bar{\mathbf{x}}) + 0\cdot x_{d_e + 1} + \cdots + 0\cdot x_D.
\end{equation*}
Then, we rotate the function $h(\mathbf{x})$ for a non-trivial constant subspace by applying a random orthogonal matrix $\mathbf{Q}$ to $\mathbf{x}$. Hence, we obtain our $D$-dimensional low-rank test function, which is given by
\begin{equation*}
    f(\mathbf{x}) = h(\mathbf{Qx}).
\end{equation*}
It is noteworthy that the first $d_e$ rows of $\mathbf{Q}$ form the basis of the effective subspace $\mathcal{T}$ of $f$, while the last $D - d_e$ rows span the constant subspace $\mathcal{T}^{\bot}.$
\begin{table}[!htb]
    \centering
    \begin{tabular}{|c|c|c|c|c|}
    \hline
    \# & Function & Effective Dimensions $d_e$ & Domain & Global Minimum \\
    \hline
    1 & low-rank Ackley \cite{Ernesto2005} & 4 & $\mathbf{x} \in [-5, 5]^4$ & 0 \\
    \hline 
    2 & low-rank Rosenbrock \cite{Surjanovic2013} & 4 & $\mathbf{x} \in [-5, 10]^4$ & 0 \\
    \hline
    3 & low-rank Shekel 5 \cite{Surjanovic2013} & 4 & $\mathbf{x} \in [0, 10]^4$ & -10.1532 \\
    \hline
    4 & low-rank Shekel 7 \cite{Surjanovic2013} & 4 & $\mathbf{x} \in [0, 10]^4$ & -10.4029 \\
    \hline
    5 & low-rank Styblinski-Tang \cite{Surjanovic2013} & 4 & $\mathbf{x} \in [-5, 5]^4$ & -156.664 \\
    \hline
    \end{tabular}
\caption{Benchmark high-dimensional low-rank test problems from \cite{cartis2020dimensionalityreductiontechniqueunconstrained, Cartis_2021}.}
\label{Low-rank Test problem}
\end{table}

\section{Additional Details}
\subsection{Sequential Domain Reduction} \label{Appendix BO_BOVAE with SDR}
To formally introduce SDR \cite{SDR}, let us denote $\mathbf{x}^{(k)} \in \mathbb{R}^D$ be the current optimal position, i.e., the centre point of the current sub-region, at iteration $k$ with each component being bounded, $x_i^{l, k} \leq x_i \leq x_i^{u, k}, i \in \{1, \ldots, D\}.$ Initially, when $k = 0$, we construct the first Region of Interest (RoI) centring at $\mathbf{x}^{(0)}$ with lower and upper bounds being
\begin{equation}\label{SDR initial bounds}
    x_i^{l, 0} = x_i^{(0)} - \frac{r_i^{(0)}}{2}, \ x_i^{u, 0} = x_i^{(0)} + \frac{r_i^{(0)}}{2}, \ i \in \{1,\ldots D\},
\end{equation}
where $r_i^{(0)}$ is the initial range value computed from the upper and lower bounds of the initial search domain. Now, suppose we are progressing from iterations $k-1$ to $k$ and that the best observations are $\mathbf{x}^{k - 1}$ and $\mathbf{x}^{k}$ up to the $(k-1)$-th and $k$-th iterations respectively. To update and contract on the RoI, we first determine an oscillation indicator along dimension $i$ at iteration $k$ as
\begin{equation}\label{oscillation indicator c}
    c_i^{(k)} = d_i^{(k)}d_i^{(k - 1)} \ \text{with} \ d_i^{(j)} = \frac{2(x_i^{(j)} - x_i^{(j - 1)})}{r^{(j - 1)}_i},  \ j = k, k - 1,
\end{equation}
where $r^{(k-1)}_i$ is the RoI size along dimension $i$ at iteration $k - 1$. Then, we normalise it as
\begin{equation}\label{normalised c}
    \hat{c}^{(k)}_i = \sqrt{|c_i^{(k)}|} \ sign(c_i^{(k)}),
\end{equation}
where $sign(\cdot)$ is the standard sign function. Then, the contraction parameter along dimension $i$ at iteration $k$ is
\begin{equation}\label{contraction parameter}
    \gamma = \frac{\gamma_p(1 + \hat{c}) + \gamma_o(1 - \hat{c})}{2},
\end{equation}
where the indices $i,k$ have been intentionally omitted for clarity and to avoid complex notations. Here, the parameter $\gamma_o$, typically set between $0.5$ and $0.7$, is a shrinkage factor to dampen oscillation. This parameter controls the reduction of the RoI, facilitating more stable and efficient convergence towards the global optimum. Meanwhiel, $\gamma_p$ indicates the pure panning behaviour and is typically set as a unity. To shrink the RoI, we utilise a zooming parameter $\eta$ to update the range along each dimension, i.e, 
\begin{equation}\label{update bounds}
    r^{(k)}_i = \lambda_i r^{(k-1)}_i, \ \text{where} \ \lambda_i = \eta + |d_i^{(k)}|(\gamma - \eta).
\end{equation}
$\lambda_i$ represents the contraction rate along dimension $i$ and $\eta$ typically lies in $[0.5, 1)$. Below, we present the Bayesian Optimisation algorithms innovatively with SDR in the ambient and the VAE-generated latent spaces.
\begin{algorithm2e}\label{BOSDR}
  \caption{Bayesian Optimisation with Sequential Domain Reduction}
  \SetAlgoLined
  \LinesNumbered
  \KwData{Initial dataset $\mathcal{D}_0 = \{\mathbf{X}_0, \mathbf{f}_0\}$, budget $B$, acquisition function $u(\cdot)$, initial search domain $\mathcal{X}$, parameters $\gamma_o, \gamma_p, \eta$, minimum region of interest size $t$, and step size $\xi$.}
  \KwResult{Minimum value $f_{min}$ found by the algorithm.}
  \textbf{Initialise} SDR by computing the initial Region of Interest (RoI) $R^{(0)}$ according to the bounds\;
  
  \For{$k = 0, \ldots, B - 1$}{
    \textbf{Fit} Gaussian process $\mathcal{GP}_k$ to current data $\mathbf{X}_k$ and $\mathbf{f}_k$\;
    \textbf{Find} the next iterate $\mathbf{x}_{k + 1} \leftarrow \argmax_{\mathbf{x}\in\mathcal{X}} u(\mathbf{x}|\mathcal{D}_k)$\;
    \textbf{Evaluate} function $f$ at $\mathbf{x}_{k+1}$, store result $f_{k+1} \leftarrow f(\mathbf{x}_{k + 1})$\;
    \textbf{Augment} the data: 
    $
      \mathbf{X}_{k + 1} \leftarrow \mathbf{X}_k \cup \{\mathbf{x}_{k + 1}\},\ \mathbf{f}_{k + 1} \leftarrow \mathbf{f}_k \cup \{f_{k+1}\}
    $
    
    \eIf{$k$ mod $\xi = 0$ and $r_i^{(k)} \geq t$}{
      \textbf{Update} RoI $R^{(k)}$ based on bounds using oscillation indicators\;
      \textbf{Trim} the updated RoI\;
    }{
      Continue to next iteration\;
    }
  }
\end{algorithm2e}
\begin{algorithm2e}[!htb]\label{BOVAE-SDR}
  \caption{BO-VAE Combined with SDR}
  \SetAlgoLined
  \LinesNumbered
  \KwData{Unlabelled dataset $\mathcal{D}_\mathbb{U} = \{\mathbf{x}_i\}_{i = 1}^M$, Initial labelled dataset $\mathcal{D}_\mathbb{L} = \{\mathbf{x}_i, f(\mathbf{x}_i)\}_{i = 1}^N$, budget $B$, initial bound $R^0$ in latent space $\mathcal{Z}$, the EI acquisition function $u(\cdot)$, an encoder-decoder pair of a VAE, $q_{\boldsymbol{\phi}}\mathbf{(z|x)}: \mathcal{X} \rightarrow \mathcal{Z}$ and $p_{\boldsymbol{\theta}}\mathbf{(x|z)}: \mathcal{Z} \rightarrow \mathcal{X}$.}
  \KwResult{Minimum value $f_{min}$ discovered by the algorithm.}
  
  \textbf{Train} the encoder $q_{\boldsymbol{\phi}}\mathbf{(z|x)}$ and decoder $p_{\boldsymbol{\theta}}\mathbf{(x|z)}$ on $\mathcal{D}_\mathbb{U}$: 
  $
    \boldsymbol{\theta}^\ast, \boldsymbol{\phi}^\ast = \argmax_{\boldsymbol{\theta}, \boldsymbol{\phi}} \mathcal{L}(\boldsymbol{\theta}, \boldsymbol{\phi}; \mathcal{D}_\mathbb{U}).
  $
  
  \textbf{Compute} the latent dataset $\mathcal{D}_{\mathbb{Z}}^0 = \{\mathbf{z}_i, f(\mathbf{x}_i)\}_{i = 1}^N$, where $\mathbf{z}_i = \mathbb{E}_{q_{\boldsymbol{\phi}^\ast}(\mathbf{z|x_i})}[\mathbf{z}]$, on $\mathcal{D}_\mathbb{L}$.

  \textbf{Initialise} SDR with the initial bound $R^0$\;
  
  \For{$k = 0, \ldots, B - 1$}{
    \textbf{Fit} GP model $h_k: \mathcal{Z} \rightarrow \mathbb{R}$ on $\mathcal{D}_\mathbb{Z}^k$\;
    
    \textbf{Solve} for the next latent point $\hat{\mathbf{z}}_k = \argmax_{\mathbf{z}} u(\mathbf{z} | \mathcal{D}_\mathbb{Z}^k)$ and reconstruct the corresponding input, $\hat{\mathbf{x}}_k \sim p_{\boldsymbol{\theta}^\ast}(\cdot|\hat{\mathbf{z}}_k)$\;
    
    \textbf{Evaluate} the objective function $f_k = f\left(\hat{\mathbf{x}}_k\right)$\;
    
    \textbf{Augment} the latent dataset $\mathcal{D}_{\mathbb{Z}}^{k + 1} \leftarrow \mathcal{D}_{\mathbb{Z}}^{k} \cup \{\hat{\mathbf{z}}_k, f_k\}$\;
    
    \textbf{Update} the search domain $R^{k + 1} \leftarrow R^k$ using the updated dataset $\mathcal{D}_{\mathbb{Z}}^{k + 1}$\;
  }
\end{algorithm2e}

\subsection{Methodology for Comparing Algorithms and Solvers}\label{Appendix Methodology for Comparing Algorithms and Solvers}
To evaluate performances of different algorithms/solvers fairly, we adopt the methodology from \cite{Cartis_2021}, using performance and data profiles as introduced in \cite{More2009}.

\textbf{Performance profiles.} A performance profile compares how well solvers perform on a problem set under a budget constraint. For a solver \(s\) and problem \(p\), the performance ratio is:
\[
r_{p,s} = \frac{M_{p,s}}{\min_{s \in \mathcal{S}} M_{p,s}},
\]
where \( M_{p,s} \) is a performance metric, typically the number of function evaluations required to meet the stopping criterion:
\[
N_p(s;\tau) = \text{\# \text{evaluations to achieve} } f^{\ast}_k \leq f^{\ast} + \tau(f^{\ast}_0 - f^{\ast}),
\]
where \(\tau \in (0,1)\) is an accuracy level. If the criterion is not met, \(N_p(s;\tau) = \infty\). The performance profile \(\pi_{s,\tau}(\alpha)\) is the fraction of problems where \(r_{p,s} \leq \alpha\), representing the cumulative distribution of performance ratios.

\textbf{Data profiles.} The data profile shows solver performance across different budgets. For a solver \(s\), accuracy level \(\tau\), and problem set \(\mathcal{P}\), it is defined as:
\[
d_{s,\tau}(\alpha) = \frac{|\{ p \in \mathcal{P}: N_p(s;\tau) \leq \alpha(n_p+1)\}|}{|\mathcal{P}|}, \ \alpha \in [0, N_g],
\]
where \(n_p\) is the problem dimension and \(N_g\) is the maximum budget. The data profile tracks the percentage of problems solved as a function of the budget.

\subsection{Soft and Hard Triplet Losses}\label{Appendix triplet loss}
In this subsection, we briefly present how \cite{grosnit2021highdimensionalbayesianoptimisationvariational} integrates the triplet loss with VAEs. We refer the reader to \cite{Hoffer2015, Wu2019} for background knowlegde for triplet deep metric loss. \cite{grosnit2021highdimensionalbayesianoptimisationvariational} introduces a parameter $\eta$ to create sets of positive \(\mathcal{D}_p(\mathbf{x}^{(b)}; \eta) = \{\mathbf{x} \in \mathcal{D}: | f(\mathbf{x}^{(b)}) - f(\mathbf{x})| < \eta\}\) and negative points \(\mathcal{D}_n(\mathbf{x}^{(b)}; \eta) = \{\mathbf{x} \in \mathcal{D}: | f(\mathbf{x}^{(b)}) - f(\mathbf{x})| \geq \eta\}\) for a base point $\mathbf{x}^{(b)}$ in a dataset $\mathcal{D}$, based on differences in function values. However, the classical triplet loss is discontinuous, which hinders GP models. To resolve this, a smooth version, the soft triplet loss, is proposed. Suppose we have a latent triplet $\mathbf{z}_{ijk} = \langle \mathbf{z}_i, \mathbf{z}_j, \mathbf{z}_k \rangle$ associated with the triplet $\mathbf{x}_{ijk} = \langle \mathbf{x}_i, \mathbf{x}_j, \mathbf{x}_k \rangle$ in the ambient space. Here, $\mathbf{z}_i$ is the latent base point. The complete expression of the soft triplet loss is \cite{grosnit2021highdimensionalbayesianoptimisationvariational}
\begin{equation*}
    \mathcal{L}_{s-trip}(\mathbf{z}_{ijk}) = \ln \left(1 + \exp(d_{\mathbf{z}}^+ - d_{\mathbf{z}}^-) \right) \omega_{ij} \omega_{ik} \times I_{\{|f(\mathbf{x}_i) - f(\mathbf{x}_j)| < \eta \ \&  \ |f(\mathbf{x}_i) - f(\mathbf{x}_k)| \geq \eta \}},
\end{equation*}
where 
\begin{equation*}
    \begin{split}
        d_{\mathbf{z}}^+ = \| \mathbf{z}_i- \mathbf{z}_j\|_p&, d_{\mathbf{z}}^- = \| \mathbf{z}_i - \mathbf{z}_k \|_p,\\
        \omega_{ij} = \frac{f_\nu\left(\eta - |f(\mathbf{x}_i) - f(\mathbf{x}_j)|\right)}{f_{\nu}(\eta)}&,
        \omega_{ik} = \frac{f_{\nu}\left(|f(\mathbf{x}_i) - f(\mathbf{x}_k)| - \eta\right)}{f_\nu(1-\eta)},
    \end{split}
\end{equation*}
for any $\mathbf{z}_j \sim q_{\boldsymbol{\phi}}(\cdot | \mathbf{x}_j), \forall \mathbf{x}_j \in \mathcal{D}_p(\mathbf{x}_i; \eta)$ and $\mathbf{z}_k \sim q_{\boldsymbol{\phi}}(\cdot | \mathbf{x}_k), \forall \mathbf{x}_k \in \mathcal{D}_n(\mathbf{x}_i;\eta)$. Here, $f_\nu (x) = \tanh \left(a/(2\nu)\right)$ is a smoothing function with $\nu$ being a hyperparameter such that $\mathcal{L}_{s-trip}(\mathbf{z}_{ijk})$ approaches the hard triplet loss since $\lim_{\nu \rightarrow 0} f_{\nu}(a) = 1$. The function $I_{\{\cdot\}}$ is a indicator function. The penalisation weights $\omega_{ij}$ and $\omega_{ik}$ are introduced to smooth out the discontinuities. Thus, the modified ELBO of a VAE trained with soft triplet loss is \cite{grosnit2021highdimensionalbayesianoptimisationvariational, Ishfaq2018}
\begin{equation*}
    \begin{split}
        \mathcal{L}_{DML} (\boldsymbol{\theta}, \boldsymbol{\phi}; \{\mathbf{x}_i, f(\mathbf{x}_i)\}_{i = 1}^N) &= \mathcal{L} (\boldsymbol{\theta}, \boldsymbol{\phi}; \{\mathbf{x}_i\}_{i = 1}^N) - \mathcal{L}_{metric} \\
        &= \sum_{n = 1}^N \left[ \mathbb{E}_{q_{\boldsymbol{\phi}}(\mathbf{z}_n | \mathbf{x}_n)} \left[\ln p_{\boldsymbol{\theta}}(\mathbf{x}_n | \mathbf{z}_n)\right] - D_{KL}\left( q_{\boldsymbol{\phi}}(\mathbf{z}_n | \mathbf{x}_n) \| p(\mathbf{z}_n)\right) \right] \\
        &\quad - \sum_{i, j, k = 1}^{N, N, N} \mathbb{E}_{q_{\boldsymbol{\phi}}(\mathbf{z}_{ijk} | \mathbf{x}_{ijk})} \left[\mathcal{L}_{s-\text{trip}}(\mathbf{z}_{ijk})\right],
    \end{split}
\end{equation*}
where $q_{\boldsymbol{\phi}}(\mathbf{z}_{ijk} | \mathbf{x}_{ijk}) = q_{\boldsymbol{\phi}}(\mathbf{z}_i | \mathbf{x}_i)q_{\boldsymbol{\phi}}(\mathbf{z}_j | \mathbf{x}_j)q_{\boldsymbol{\phi}}(\mathbf{z}_k | \mathbf{x}_k).$

The BO-VAE algorithm with the soft triplet loss as the chosen deep metric loss is outlined \ref{BOVAE with DML}. We note that Algorithm \ref{BOVAE with DML} is not implemented with SDR in the latent space, as experiments have shown that SDR and DML methods conflict with each other in excluding the global optimum. Addressing this conflict when implementing SDR in DML-structured latent spaces is left as future work.
\begin{algorithm2e}[!htb]\label{BOVAE with DML}
  \caption{Retraining BO-VAE Algorithm with DML}
  \SetAlgoLined
  \LinesNumbered
  \KwData{Labelled dataset $\mathcal{D}_\mathbb{L}^{l = 1} = \{\mathbf{x}_i, f(\mathbf{x}_i)\}_{i = 1}^N$, unlabelled dataset $\mathcal{D}_\mathbb{U} = \{\mathbf{x}_i\}_{i = 1}^M$, budget $B$, periodic frequency $q$, EI acquisition function $u(\cdot)$, the encoder and decoder models from a VAE, $q_{\boldsymbol{\phi}}\mathbf{(z|x)}: \mathcal{X} \rightarrow \mathcal{Z}$ and $p_{\boldsymbol{\theta}}\mathbf{(x|z)}: \mathcal{Z} \rightarrow \mathcal{X}$.}
  \KwResult{Minimum function value $f_{min}$ found by the algorithm.}

  \textbf{Pre-train} the VAE model $V_{\mathcal{D}_\mathbb{L}}^{l=0}$ with $\mathcal{D}_\mathbb{U}$: $
  \boldsymbol{\theta}^\ast_0, \boldsymbol{\phi}^\ast_0 = \argmax_{\boldsymbol{\theta}, \boldsymbol{\phi}} \mathcal{L}(\boldsymbol{\theta}, \boldsymbol{\phi}; \mathcal{D}_\mathbb{U})
  $

  Set $\boldsymbol{\theta}^\ast_1 \leftarrow  \boldsymbol{\theta}^\ast_0$, $\boldsymbol{\phi}^\ast_1 \leftarrow \boldsymbol{\phi}^\ast_0$, $V_{\mathcal{D}_\mathbb{L}}^{l=1} \leftarrow V_{\mathcal{D}_\mathbb{L}}^{l=0}$\;

  \For{$l = 1$ \text{to} $L \equiv \lceil B/q \rceil$}{
    \textbf{Train} the VAE model $V_{\mathcal{D}_\mathbb{L}^l}^{l}$ on $\mathcal{D}_\mathbb{L}$:
    $
    \boldsymbol{\theta}^\ast_l, \boldsymbol{\phi}^\ast_l = \argmax_{\boldsymbol{\theta}, \boldsymbol{\phi}} \mathcal{L}_{DML}(\boldsymbol{\theta}, \boldsymbol{\phi}; \mathcal{D}_\mathbb{L}^l)
    $

    \textbf{Compute} the latent dataset 
    $
    \mathcal{D}_\mathbb{Z}^l = \{ \mathbf{z}_i, f(\mathbf{x}_i) \}^{N + l \cdot q}_{i = 1} = \{ \mathbb{E}_{q_{\boldsymbol{\phi}^\ast_l}(\mathbf{z|x_i})}[\mathbf{z}], f(\mathbf{x}_i) \}^{N + l \cdot q}_{i = 1}
    $
    
    \textbf{Initialise} $\mathcal{D}_\mathbb{L}^{l;k=0} \leftarrow \mathcal{D}_\mathbb{L}^{l}$ and $\mathcal{D}_\mathbb{Z}^{l;k=0} \leftarrow \mathcal{D}_\mathbb{Z}^{l}$\;

    \For{$k = 0$ \text{to} $q-1$}{
      \textbf{Fit} a Gaussian Process (GP) model $h_{l;k}: \mathcal{Z} \rightarrow \mathbb{R}$ on
      $
      \mathcal{D}_\mathbb{Z}^{l;k} = \{ \mathbf{z}_i, f(\mathbf{x}_i) \}^{N + l \cdot q + k}_{i = 0}
      $
      
      \textbf{Solve} for the next latent point:
      $
      \hat{\mathbf{z}}_{l; k + 1} = \argmax_{\mathbf{z}} u(\mathbf{z}|\mathcal{D}_\mathbb{Z}^{l;k})
      $
      
      \textbf{Obtain} the new sample $\hat{\mathbf{x}}_{l; k + 1}$: 
      $
      \hat{\mathbf{x}}_{l; k + 1} \sim p_{\boldsymbol{\theta}^\ast_l}(\cdot| \hat{\mathbf{z}}_{l; k + 1})
      $

      \textbf{Evaluate} the objective function at the new sample: $f(\hat{\mathbf{x}}_{l; k + 1})$\;

      \textbf{Augment} the datasets:
      \[
      \mathcal{D}_\mathbb{L}^{l; k+1 } \leftarrow \mathcal{D}_\mathbb{L}^{l;k} \cup \{ \hat{\mathbf{x}}_{l; k+ 1}, f(\hat{\mathbf{x}}_{l; k+1}) \}
      ,
      \mathcal{D}_\mathbb{Z}^{l; k+1} \leftarrow \mathcal{D}_\mathbb{Z}^{l;k} \cup \{ \hat{\mathbf{z}}_{l; k + 1}, f(\hat{\mathbf{x}}_{l; k + 1}) \}
      \]
    }

    \textbf{Augment} the outer-loop datasets:
    $
    \mathcal{D}_\mathbb{L}^{l+1} \leftarrow \mathcal{D}_\mathbb{L}^{l;q}, \mathcal{D}_\mathbb{Z}^{l+1} \leftarrow \mathcal{D}_\mathbb{Z}^{l;q}
    $
  }
\end{algorithm2e}

\section{Additional Details for Section 4 Numerical Experiments} \label{Appendix Additional Details for Section 4 Numerical Experiments}
To train VAEs more robustly, we introduce an additional weight $\beta$ to the VAE ELBO, known as beta-annealing approach \cite{bowman2016generatingsentencescontinuousspace, burgess2018understandingdisentanglingbetavae}. The modified ELBO becomes
\begin{equation*}
    \begin{split}
       \mathcal{L}(\mathbf{\boldsymbol{\theta}, \boldsymbol{\phi}}; \mathbf{x}) & =   \ln p_{\boldsymbol{\theta}}(\mathbf{x}) - \beta D_{KL} [q_{\boldsymbol{\phi}}(\mathbf{z|x}) \| p_{\boldsymbol{\theta}}(\mathbf{z|x})] \\
       & = \mathbb{E}_{q_{\boldsymbol{\phi}}(\mathbf{z|x})} [\ln p_{\boldsymbol{\theta}}(\mathbf{x|z})] - \beta D_{KL} [q_{\boldsymbol{\phi}}(\mathbf{z|x}) \| p(\mathbf{z})],
    \end{split}
\end{equation*}
where $\beta \geq 0.$ The use of $\beta$ is a trade-off between reconstruction accuracy and the regularity of the latent space and to avoid the case of the vanishing KLD term, where no useful information is learned \cite{bowman2016generatingsentencescontinuousspace, fu2019cyclicalannealingschedulesimple}. A common approach to implementing the beta-annealing technique \cite{bowman2016generatingsentencescontinuousspace} involves initialising $\beta$ at $0$ and gradually increasing it in uniform increments over equal intervals until $\beta$ reaches $1.$ 
\paragraph{Experimental Configurations for Test Set 1} The basic training details of the VAE used in the experiments are given in Table \ref{VAE Training Table}.
\begin{table}[!htb]
    \centering
    \begin{tabular}{ | c | c| c |  c | c | c |}
    \hline 
     Epochs & Optimiser & Learning Rate & Batch Size & $\left(\beta_i, \beta_f, \beta_s, \beta_a\right)$ & $M$  \\
    \hline
    $300$ & Adam & $1 \times 10^{-3}$ & $1024$ & $(0, 1, 10, 0.1)$ & $50000$\\
    \hline
    \end{tabular}
    
\caption{\small The (pre-)training details of the VAE used for Test Set 1. $\beta_i$ and $\beta_f$ are the initial and final values of $\beta$ for $\beta-$VAEs respectively. The annealing approach is to increase the weight $\beta$ by $\beta_a$ every $\beta_s$ epochs. $M$ is the size of \( \mathcal{D}_{\mathbb{U}} \).}
\label{VAE Training Table}
\end{table}

We highlight two ingredients in the implementation of the algorithms. 
\begin{enumerate}
    \item The first thing involves the VAE pre-training. The models are pre-trained according to the details in Table \ref{VAE Training Table}. It is crucial that training samples are drawn with high correlations to construct the VAE training dataset. For instance, samples can be generated from a multivariate normal distribution with a large covariance matrix. This approach facilitates the VAE in learning a meaningful low-dimensional data representation.
    
    \item  The second one involves constructing the latent datasets for a sample-efficient BO procedure, as it would be computationally inefficient to use the entire VAE training dataset. Therefore, instead of using the entire \( \mathcal{D}_{\mathbb{L}} \), we utilise only 1\% of it by uniformly and randomly selecting \( N \) points, where \( N \) represents 1\% of the size of \( \mathcal{D}_{\mathbb{U}} \) at the current retraining stage \( l \). 
\end{enumerate}

The SDR setting is: $\gamma_o = 0.7$, $\gamma_p = 1.0$, $\eta = 0.9$, $t = 0.5$, $\xi = 1$. The initial search domain $R^0$ for Algorithms \ref{general BOVAE} and \ref{BOVAE-SDR} is $[-5, 5]^d$. For Algorithm \ref{BOVAE with DML}, the hyperparameters $\eta$ and $\nu$ are set to be $0.01$ and $0.2$ respectively. For the retraining stage, we use Table \ref{Chap 3 Retrain VAE details Table} as the common setup.
\begin{table}[!htb]
    \centering
    \begin{tabular}{| c | c | c | c |  c |}
    \hline 
     Epochs & Optimiser & Learning Rate & Batch Size & Beta-annealing \\
    \hline
     $2$ & Adam & $1 \times 10^{-3}$ & $256$ & No  \\
    \hline
    \end{tabular}
    
\caption{The retraining details of the VAE used for Test Set 1.}
\label{Chap 3 Retrain VAE details Table}
\end{table}

\paragraph{Experiment Configurations for Test Set 2}
The pre-training and retraining details for this VAE are consistent with as before, as shown in Table \ref{VAE Training Table} and Table \ref{Chap 3 Retrain VAE details Table}, respectively. In addition to the two key implementation details for BO-VAE algorithms listed in Appendix \ref{Appendix Additional Details for Section 4 Numerical Experiments}, it is important to note that the test problem domains must be scaled to \([-1, 1]^D\) for a fair comparison with REMBO. This adjustment is due to the domain scaling used in constructing Test Set 2. The specific experimental configurations for each BO-VAE algorithm are consistent with as before.

\end{document}